\documentclass{article}
\usepackage{amssymb}
\usepackage{amsmath}
\usepackage{latexsym}
\usepackage{mathrsfs}
\usepackage{graphicx}
\usepackage[all]{xy}
\usepackage[CJKbookmarks=true]{hyperref}
\newtheorem{Th}{Theorem}
\newtheorem{Lemma}[Th]{Lemma}
\newtheorem{Coro}[Th]{Corollary}
\title{Comodule Structures, Equivariant Hopf Structures, and Generalized Schubert Polynomials}
\author{XIONG Rui}
\begin{document}

\def\CH{\operatorname{CH}}
\maketitle

\begin{abstract}
  In this article, the comodule structure of Chow rings of Flag manifolds $\CH(G/B)$ is described by Schubert cells.
  Its equivariant version gives rise to a Hopf structure of the equivariant cohomology of flag manifolds $H_B^*(G/B)$.
  We get two identities of generalized Schubert polynomials as explanations of the geometric facts.
\end{abstract}

I would politely express my gratitude to Victor Petrov, Haibao Duan, Xuezhi Zhao, Paul Zinn-Justin and Shintar\^o Kuroki for discussion.

\section{Main Results}

\paragraph{Comodule structure of $\CH(G/B)$. }
Let $G$ be a reductive group, $B$ its Borel subgroup.
Denote the Weyl group by $W=N_G(T)/T$, the length function $\ell$.
For $u,v,w\in W$, we write $w=u\odot v$ if $w=uv$ and $\ell(w)=\ell(u)+\ell(v)$.
It follows from the Bruhat decomposition (for example, \cite{SpringerLinear})
that the flag $G/B$ has a cellular structure by $B$-orbits, $\{BwB/B:w\in W\}$.
So the Chow ring $\CH(G/B)$ is freely generated by the fundamental class of $\{\overline{BwB/B}:w\in W\}$.
We denote $[\Sigma_w]$ the class corresponding to $Bw_0wB/B$, so that $[\Sigma_w]\in \CH^{2\ell(w)}(G/B)$.
Denote $\pi:G\to G/B$ the natural projection.


\begin{Th}\label{ComoduleTh}The map induced by the left group action $G\times G/B\to G/B$ is given by
$$\CH(G/B)\longrightarrow \CH(G)\otimes \CH(G/B)\qquad [\Sigma_w]\longmapsto \sum_{w=u\odot v} \pi^*[\Sigma_u]\otimes [\Sigma_v]. $$
\end{Th}

It can be generalized to any parabolic subgroup $P$.
Assume $P$ correspondent to the subset $\Theta$ of simple reflection generators of $W$.
Now the $B$-orbit $\{BwP/P:w\in W^P\}$ forms a cellular structure,
where $W^P$ the set of $w\in W$ with $\ell(ws)=\ell(w)+1$ for all $s\in \Theta$.
We similarly denote $[\Sigma_w]=[\overline{Bw_0wP/P}]$.
Since $G/B\to G/P$ is cellular and $\CH(G/P)\to \CH(G/B)$ is injective, we can conclude the following corollary.

\begin{Coro}The map induced by the left group action $G\times G/P\to G/P$ is given by the same formula as in theorem \ref{ComoduleTh}.
\end{Coro}

It follows from Grothendieck \cite{grothendieck1958torsion} that $\pi$ induces a surjection onto $\CH(G)$.
Hence as a partial corollary, we have the following description of Hopf structure of Chow ring.

\begin{Th}\label{HopfStructureTh}The Hopf structure of $\CH(G)$ can be described by
$$\Delta:\CH(G)\longrightarrow \CH(G)\otimes \CH(G)\qquad
\pi^*[\Sigma_w]\longmapsto \sum_{w=u\odot v} \pi^*[\Sigma_u] \otimes \pi^*[\Sigma_v], $$
and the antipode
$$S: \CH(G)\longrightarrow \CH(G)\qquad \pi^*[\Sigma_w]\longmapsto (-1)^{\ell(w)}\pi^*[\Sigma_{w^{-1}}]. $$
\end{Th}

From the paper of Kac \cite{kac1985torsion},
we know that image of $\pi^*$ is far from being injection, but this formula provides a uniform description.
Historically, the determination of Hopf structure of compact Lie groups are quite tough.
The classic approach is the cohomology operators, see for example \cite{ishitoya1976hopf}.
The determination based on computer-assisted computation can be found in Duan and Zhao \cite{duan2014schubert}, \cite{duan2015schubert}.
The method using motives can be found in Petrov and Semenov \cite{petrov2019hopftheoretic}.


\paragraph{Equivariant version. }
\def\hq{\backslash\!\backslash}
\def\pt{\mathsf{pt}}
In general, for a topological group $G$ of CW-type, $X$ a $G$-space,
we denote $X_G=EG\times_G X$, where $EG$ is the total space of Milnor's universal $G$-bundle.
Note that if $X$ is acted freely by $G$, then $X_G$ is homotopy equivalent to the orbit space $G\backslash X$.
For a left $G$-set $X$ and a right $G$-set $Y$, we denote
$$X\underset{G}*Y=(X\times Y)_G=EG\times_G(X\times Y). $$
Under this notation, $X_G=\pt \underset{G}* X$.
We denote $H^*_G(X)=H^*(X_G)$ and call it the $G$-equivariant cohomology.

Now turn to the case of reductive group. We can consider the map induced by group multiplication and inverse
$$\mu:\pt \underset{B}*G\underset{B}*G\underset{B}* \pt\to \pt \underset{B}* G\underset{B}* \pt,\qquad
\nu: \pt \underset{B}* G\underset{B}* \pt \longrightarrow \pt \underset{B}* G\underset{B}* \pt. $$
Note that $\pt \underset{B}* G\underset{B}* \pt= (G/B)_B$.
Denote the $[\Sigma_w]_B$ the equivariant fundamental class in $H^*_B(G/B)=H^*(\pt \underset{B}* G\underset{B}* \pt)$
of the $B$-equivariant subvariety $\overline{w_0Bw_0wB/B}$.
Consider the following two maps
$$\begin{array}{c}
\pi_1: \pt \underset{B}*G\underset{B}*G\underset{B}* \pt\to \pt \underset{B}*G\underset{B}*G\underset{G}* \pt
= \pt \underset{B}*G\underset{B}* \pt,\\
\pi_2: \pt \underset{B}*G\underset{B}*G\underset{B}* \pt\to \pt \underset{G}*G\underset{B}*G\underset{B}* \pt
= \pt \underset{B}*G\underset{B}* \pt.\\
\end{array}$$


\begin{Th}\label{EqHopfStructureTh}The map induced by $\mu$ is given by
$$\mu^*:H^*(\pt \underset{B}* G\underset{B}* \pt)\longrightarrow H^*_B(\pt \underset{B}*G\underset{B}*G\underset{B}* \pt)\qquad
[\Sigma_w]_B\longmapsto \sum_{w=u\odot v}\pi_1^*[\Sigma_u]_B\smile \pi_2^*[\Sigma_v]_B, $$
and the map induced by $\nu$ is given by
$$\nu^*:H^*(\pt \underset{B}* G\underset{B}* \pt)\longrightarrow H^*(\pt \underset{B}* G\underset{B}* \pt)\qquad
[\Sigma_w]_B\longmapsto  (-1)^{\ell(w)}[\Sigma_{w^{-1}}]_B. $$
\end{Th}

\paragraph{Generalized Schubert Polynomials. }
Fix a choice of $T\cong (\mathbb{C}^\times)^n$, it is proved by Borel in \cite{borel1961sous} that
$$H_B^*(G/B;\mathbb{Q})=\frac{\mathbb{Q}[x_1,\ldots,x_n,t_1,\ldots,t_n]}{\left<f(x)-f(t):f(X)\in \mathbb{Q}[X_1,\ldots,X_n]^W\right>}. $$
Then $[\Sigma_w]_B$ is uniquely presented by a polynomial $\mathfrak{S}_w(x,t)$ module the ideal factored as above.
We fix a choice of the polynomials, and call them \emph{generalized Schubert polynomials}.
Note that we follow the classic convention, writing the $t$ the indeterminant from $H_B^*(\pt)$.
As a result, $\mathfrak{S}_w(x)=\mathfrak{S}_w(x,0)$ is the classic Kostant polynomial \cite{billey1999kostant}, but
up to some element in the ideal generated by polynomials of $\mathbb{Q}[X_1,\ldots,X_n]^W$ with zero constant term.

\begin{Th}\label{doubleidentities}On the generalized Schubert polynomials, we have
$$\mathfrak{S}_w(x,t)
\equiv\sum_{w=u\odot v}\mathfrak{S}_v(x,y)\mathfrak{S}_{u}(y,t) \hfill \mod
\left<{{f(x)-f(y)}\atop{f(t)-f(y)}}:f\in \mathbb{Q}[X]^W\right>$$
and
$$\mathfrak{S}_w(x,t)\equiv (-1)^{\ell(w)} \mathfrak{S}_{w^{-1}}(t,x) \mod \left<f(x)-f(t):f\in \mathbb{Q}[X]^W\right>.$$
\end{Th}

\begin{Coro}\label{doublecomputation}In particular, by taking $y=0$,
$$\mathfrak{S}_w(x,t)
\equiv\sum_{w=u\odot v}(-1)^{\ell(u)}\mathfrak{S}_v(x)\mathfrak{S}_{u^{-1}}(t) \hfill \mod
\left<\mathbb{Q}[x]^W_+\right>+\left<\mathbb{Q}[t]^W_+\right>$$
where $\mathbb{Q}[t]^W_+$ stands the set of polynomial of $\mathbb{Q}[t]^W$ without constant term.
\end{Coro}

\def\GL{\operatorname{GL}}

In the case $G=\GL_n$, we can find a stable choice for each $\mathfrak{S}_{w}(x,t)$,
the classic double Schubert polynomial \cite{fomin1994schubert},
i.e. a polynomial with each monomial in $x$ no great than $x_1^{n-1}\cdots x_{n-1}$ in dominant order.
In this case $\mathfrak{S}_w(x)=\mathfrak{S}_w(x,0)$ is exactly the Schubert polynomials defined by
Lascoux and Sch{\"u}tzenberger \cite{lascoux1982structure}.
This choice is called stable since we can pass to the infinite flag manifold, in which case, the ideal tends to zero.
So the above formula recovers the three identities of classic Schubert polynomials which were proved by combinatorial method only,
see \cite{fomin1996yang} and \cite{lam2018stable} (note that the sign convention of $t$ in old papers is distinct from ours).

For other classic groups, Kirillov \cite{kirillov2015double} proved that there is a choice which makes Corollary \ref{doublecomputation}
holds without passing to the ideal.

\section{Geometric Part}

\paragraph{Reduction to equivariant case. }
Firstly, $\CH(G/B)=H^*(G/B)$ due to the cellular structure, and as we remarked, the image of $H^*(G/B)\to H^*(G)$ is $\CH(G)$,
so we can use the topological argument.
Let us consider the diagram
$$\xymatrix{
\pt \times G\times G\underset{B}* \pt \ar[r]\ar[d] & \pt \times  G\underset{B}* \pt\ar[d]\\
\pt \underset{B}*G\underset{B}*G\underset{B}* \pt \ar[r]_{\mu} & \pt \underset{B}*G\underset{B}* \pt
}\qquad
\xymatrix{
\pt \times G\times \pt \ar[r]\ar[d]& \pt \times G \times \pt\ar[d]\\
\pt \underset{B}* G\underset{B}* \pt \ar[r]_{\nu}& \pt \underset{B}* G\underset{B}* \pt
}$$
So if we know Theorem \ref{EqHopfStructureTh} is true, then to prove Theorem \ref{ComoduleTh} and Theorem \ref{HopfStructureTh},
it suffices to prove
$$H^*(\pt \underset{B}* G\underset{B}* \pt)\longrightarrow H^*(G)\qquad
[\Sigma_w]_B\longmapsto \pi^*[\Sigma_w], $$
and
$$H^*(\pt \underset{B}* G\underset{B}*G\underset{B}* \pt)\longrightarrow H^*(G\times G/B)\qquad
\begin{cases}
\pi_1^*[\Sigma_u]_B\longmapsto \pi^*[\Sigma_u]\\
\pi_2^*[\Sigma_v]_B\longmapsto [\Sigma_v].
\end{cases} $$
For the first map is clear, since it factors through
$$H^*(\pt \underset{B}* G\underset{B}* \pt)=H_B^*(G/B)\stackrel{\epsilon}\longrightarrow H^*(G/B)\stackrel{\pi^*}\longrightarrow H^*(G), $$
where $\epsilon$ is the augment map which maps $[\Sigma_w]_B$ to $[\Sigma_w]$.
For the second, note that
$$\xymatrix{
\pt \times G\times G\underset{B}* \pt \ar[r]&
\pt \times G\underset{B}* G\underset{B}* \pt\ar[r]^{}\ar[d]
&\pt \times G\underset{B}* G\underset{G}* \pt \ar[d]&G/B\ar[d]\ar[l]_{\sim}\\
& \pt \underset{B}*G\underset{B}*G\underset{B}* \pt\ar[r]^{\pi_1}
&\pt \underset{B}*G\underset{B}*G\underset{G}* \pt& (G/B)_B\ar[l]_{\sim}}$$
$$\xymatrix{
\pt \times G\times G\underset{B}* \pt \ar[r]&
\pt \underset{B}*G\times G\underset{B}* \pt\ar[r]^{}\ar[d]
&\pt \underset{G}*G\times G\underset{B}* \pt \ar[d]&G/B\ar[d]\ar[l]_{\sim}\\
& \pt \underset{B}*G\underset{B}*G\underset{B}* \pt\ar[r]^{\pi_2}
&\pt \underset{G}*G\underset{B}*G\underset{B}* \pt& (G/B)_B\ar[l]_{\sim}}$$
The first row of each diagram is exactly the first (second respectively) projection.

\paragraph{Reduction to polynomials. }
Note $B$ is homotopy equivalent to $T$, so we can freely exchange $B$ by $T$.
It is well-known that $H^*(BT)$ is naturally isomorphism to the ring of polynomials in characters of $T$.
If we fix a choice of $T\cong (\mathbb{C}^\times)^n$, then $H^*(BT)=\mathbb{Z}[t_1,\ldots,t_n]$.
More exactly, $t_i$ is the $i$-th projection of $(\mathbb{C}^\times)^n$.
Besides, it is known that $H^*(BG;\mathbb{Q})=H^*(BT;\mathbb{Q})^W=\mathbb{Z}[t_1,\ldots,t_r]^W$ the ring fixed by $W$,
for example \cite{hsiang2012cohomology}.

Since we will use the details of the argument of Borel \cite{borel1961sous}, let me firstly briefly repeat them.
As $BT$ can be picked to be $EG/T$, so $BT\to BG$ is a fibre bundle of fibre $G/T$.
Since the $BG$ and $G/T$ are all only of even dimensions, so the Serre--Leray spectral sequence collapses at the second page.
In particular, as $H^*(BG;\mathbb{Q})$-module
$$H^*(BT;\mathbb{Q})=H^*(BG;\mathbb{Q})\otimes H^*(G/T;\mathbb{Q}). $$
So
$$\begin{array}{rl}
H^*(G/T;\mathbb{Q})&=H^*(BT;\mathbb{Q})\otimes_{H^*(BG;\mathbb{Q})} \mathbb{Q}\\
& =H^*(BT;\mathbb{Q})/\left<H^+(BG;\mathbb{Q})\right>=\dfrac{\mathbb{Q}[x_1,\ldots,x_n]}{\left<\mathbb{Q}[x]^W_+\right>}.
\end{array}$$
Actually, $H^*(BT;\mathbb{Q})$ is free as $H^*(BG;\mathbb{Q})$-module of rank $\dim H^*(G/T;\mathbb{Q})=|W|$.
By arguing the same but with fibre at each point,
we see that any choice of basis of $H^*(BT;\mathbb{Q})$ over $H^*(BG;\mathbb{Q})$ restricting at each fibre to be a $\mathbb{Q}$-basis.

Then, consider
$$\xymatrix{
&\pt\underset{T}* G\underset{T}* \pt\ar[r]\ar[d] & \pt\underset{G}* G\underset{T}* \pt \ar[d]&BT\ar[l]_{\simeq}\\
BT\ar[r]_{\simeq}&\pt\underset{T}* G\underset{G}* \pt\ar[r] & \pt\underset{G}* G\underset{G}* \pt&BG.\ar[l]_{\simeq}}$$
Then the Harish--Leray theorem shows that
$$\begin{array}{rl}
H^*_T(G/T;\mathbb{Q})&=H^*(\pt\underset{T}* G\underset{T}* \pt;\mathbb{Q})\\
&=H^*(BT;\mathbb{Q})\otimes_{H^*(BG;\mathbb{Q})} H^*(BT;\mathbb{Q})\\
&=\dfrac{\mathbb{Q}[x_1,\ldots,x_n,t_1,\ldots,t_n]}{\left<f(x)-f(t): f\in \mathbb{Q}[X]^W\right>}.
\end{array}$$

Lastly, learning from the Borel's augment above, we should consider
$$\xymatrix{
&\pt\underset{T}* G\underset{T}* G\underset{T}* \pt\ar[r]\ar[d] & \pt\underset{G}* G\underset{T}* G\underset{T}* \pt \ar[d]&
\pt\underset{T}* G\underset{T}* \pt\ar[l]_{\sim}\\
BT\ar[r]_{\simeq}&\pt\underset{T}* G\underset{G}* G\underset{G}* \pt\ar[r] & \pt\underset{G}* G\underset{G}* G\underset{G}* \pt&BG.\ar[l]_{\simeq}}$$
Then still the Harish--Leray theorem shows that
$$\begin{array}{rl}
H^*_T(G\times_TG/T; \mathbb{Q})&=H^*(\pt\underset{T}* G\underset{T}* G\underset{T}* \pt;\mathbb{Q})\\
&=H^*(BT;\mathbb{Q})\otimes_{H^*(BG;\mathbb{Q})} H^*(\pt\underset{T}* G\underset{T}* \pt;\mathbb{Q})\\
&=H^*(BT;\mathbb{Q})\otimes_{H^*(BG;\mathbb{Q})} H^*(BT;\mathbb{Q})\otimes_{H^*(BG;\mathbb{Q})} H^*(BT;\mathbb{Q})\\
&=\dfrac{\mathbb{Q}[x_1,\ldots,x_n,y_1,\ldots,y_n,t_1,\ldots,t_n]}{\left<{{f(x)-f(y)},{f(t)-f(y)}}:f\in \mathbb{Q}[X]^W\right>}.
\end{array}$$
To be exact, let us fix the choice of where the indeterminants $\{t_i,y_i,x_i\}$ from,
$$\begin{array}{r@{}l}
H^*(\pt\underset{T}* G\underset{G}* G\underset{G}* \pt)=H^*(BT)=&\mathbb{Z}[t_1,\ldots,t_n],\\
H^*(\pt\underset{G}* G\underset{T}* G\underset{G}* \pt)=H^*(BT)=&\mathbb{Z}[y_1,\ldots,y_n],\\
H^*(\pt\underset{G}* G\underset{G}* G\underset{T}* \pt)=H^*(BT)=&\mathbb{Z}[x_1,\ldots,x_n].
\end{array}$$

\begin{Th}The $H^*_T(G/T;\mathbb{Q})=H^*(BT;\mathbb{Q})\otimes_{H^*(BG;\mathbb{Q})} H^*(BT;\mathbb{Q})$
forms a Hopf algebra over $H^*(BT;\mathbb{Q})$ under product the cup product $\smile$,
coproduct $\mu^*$ with antipole $\nu^*$.
\end{Th}

Now, consider the map
$$\pt\underset{T}* G\underset{T}* G\underset{T}* \pt\longrightarrow  \pt\underset{T}* G\underset{G}* G\underset{T}* \pt
=\pt\underset{T}* G\underset{T}* \pt. $$
This is nothing but our $\mu$.
Since the order of exchanging $\underset{T}*$ by $\underset{G}*$ does not matter, so
in conclusion, under our computation,
$$\mu^*:H^*(\pt \underset{B}* G\underset{B}* \pt)\longrightarrow H^*_B(\pt \underset{B}*G\underset{B}*G\underset{B}* \pt)\qquad
f(x,t)\longmapsto f(x,t). $$
Furthermore,
$$\begin{array}{rl}
\pi_1^*:&H^*(\pt \underset{B}* G\underset{B}* \pt)\longrightarrow H^*_B(\pt \underset{B}*G\underset{B}*G\underset{B}* \pt)\qquad
f(x,t)\longmapsto f(y,t),\\
\pi_2^*:&H^*(\pt \underset{B}* G\underset{B}* \pt)\longrightarrow H^*_B(\pt \underset{B}*G\underset{B}*G\underset{B}* \pt)\qquad
f(x,t)\longmapsto f(x,y).
\end{array}$$

To describe $\nu$, it is not very hard. Note that
$$\xymatrix{
\pt \underset{T}* G\underset{T}* \pt \ar[r]\ar[d]_{\nu} & \pt \underset{T}* G\underset{G}* \pt\ar[d]& \pt \underset{T}* \pt \ar[d]^{B\nu}\ar[l]_{\simeq}\\
\pt \underset{T}* G\underset{T}* \pt \ar[r] & \pt \underset{G}* G\underset{T}* \pt& \pt \underset{T}* \pt \ar[l]_{\simeq}\\
}$$
where $B\nu$ is the map $BT\to BT$ induced by the inverse $T\to T$. As a result,
$$\nu^*:H^*(\pt \underset{B}* G\underset{B}* \pt)\longrightarrow H^*(\pt \underset{B}* G\underset{B}* \pt)
\qquad f(x,t)\longmapsto f(-t,-x). $$

By an easy spectral sequence argument, we see $H^*(\pt \underset{B}* G\underset{B}* \pt)$ is free abelian.
So to check the expression of $\mu^*$ and $\nu^*$, it is harmless to reduce to over $\mathbb{Q}$.
As a result, Theorem \ref{EqHopfStructureTh} follows from Theorem \ref{doubleidentities}.

\section{Algebraic Part}

\paragraph{Demazure operators and Localization. }
Let us denote the longest element by $w_0\in W$, and simply write $(-1)^w$ for $(-1)^{\ell(w)}$.
Denote $\Sigma=\{\alpha_i:i\in I\}$ the set of simple root, and $s_i: i\in I$ the corresponding simple reflection of $\alpha_i$.
We define the \emph{Demazure operator} over $H^*(BT)$ to be
$$\partial_i: H^*(BT)\longrightarrow H^{*-2}(BT)\qquad  f\longmapsto\frac{f-s_if}{\alpha_i}. $$
Note that $\partial_i^2=0$ for each $i$.
It is classic that $\{\partial_i\}$ satisfies the Braid relations,
so for any $w=s_{i(1)}\ldots s_{i(r)}\in W$ a reduced decomposition, the operator
$$\partial_w=\partial_{i(1)}\circ \ldots \circ \partial_{i(n)}$$
does not depend on the choice of reduced decomposition, see for example \cite{hiller1982geometry}.
We also call $\partial_v$ the \emph{Demazure operator}.

Since we work in more variables, we will use $\partial_v^x$ to stand the operator with respect to $x$.
It is known that $\partial_v^x$ can be reduced to $H_B^*(G/B)$, and well-behaved on Schubert cells,
$$\partial_v^x: H^*_B(G/B)\longrightarrow H^{*-2\ell(v)}(G/B)\qquad
[\Sigma_w]_B\longmapsto \begin{cases}
[\Sigma_{wv^{-1}}]_B, & \ell(wv^{-1})+\ell(v)=\ell(w), \\
0, & \text{otherwise}.
\end{cases}$$
In other word,
$$\partial_v^x\mathfrak{S}_w(x,t)\equiv \begin{cases}
\mathfrak{S}_{wv^{-1}}(x,t), & \ell(wv^{-1})+\ell(v)=\ell(w), \\
0, & \text{otherwise}.
\end{cases} \mod
\left<{{f(x)-f(y)}\atop{f(t)-f(y)}}:f\in \mathbb{Q}[X]^W\right>. $$
The nonequivariant form is proved in \cite{bernstein1982schubert},
and the equivariant form and further discussion can be found in \cite{kaji2015presentations}.

Consider the localization map
$$\cdot |_{w}: H^*_T(G/B)\to H^*_T(wB/B)=H^*_T(\pt). $$
By the choice of $[\Sigma_w]_B$, for $u\in W$,
$$\mathfrak{S}_w(x,ux)\neq 0\Longrightarrow u\geq w. $$

Actually, the above two properties in above two paragraphs as well as $\mathfrak{S}_{1}(x,t)=1$ characterize $\mathfrak{S}_{w}(x,t)$.

\def\NH{\operatorname{NH}}
\paragraph{Affine Nil-Hecke Algebra. }
Denote the affine Nil-Hecke Algebra
$\NH_W$ the algebra generated in $H^*(BT)$ by the Demazure operators $\{\partial_w:w\in W\}$ and left multiplication of element in $H^*(BT)$.
Abuse of notation, we write $f(x)$ by multiplication by $f(x)$.
Since
$$\partial_i (fg)=\frac{fg-s_if s_ig}{\alpha}=\frac{fg-s_ifg}{\alpha}+\frac{s_ifg-s_ifs_ig}{\alpha}
=(\partial_i f)g+(s_if)(\partial_i g), $$
we have the following Leibniz in $\NH_W$,
$$\partial_i f=(\partial_if)+(s_if)\partial_i. $$

\begin{Lemma}\label{TotalLeibniz}In $\NH_W$, we have the following identity,
$$(-1)^{w_0}\partial_{w_0} F(w_0x)=\sum_{w\in W} \big(\partial_{ww_0}F(x)\big)\cdot (-1)^{w}\partial_{w}.$$
Or, equivalently, in term of polynomials,
$$(-1)^{w_0}\partial_{w_0} \big(F(w_0x)G(x)\big)=\sum_{w\in W} \partial_{ww_0}F(x)\cdot (-1)^{w}\partial_{w}G(x).$$
\end{Lemma}

A geometric proof of this theorem will be given in the next section.
But here I present the pure algebraic proof.

\textit{Proof. }
We can assume $(-1)^{w_0}\partial_{w_0} F(w_0x)=\sum_{w\in W} c_{w}(x)\cdot (-1)^{w}\partial_{w}$.
Since whenever we apply $\partial_i$, the left hand side will vanish, we have
$$\sum_{w\in W} \big(\partial_i c_{w}(x)\big)\cdot (-1)^{w}\partial_{w}+\sum_{w\in W} c_{w}(s_ix)\cdot (-1)^{w}\partial_i\partial_{w}=0. $$
In other word, $\partial_i c_{w}(x)=\begin{cases}
c_{s_iw}(s_ix), & s_iw<w,\\
0,& s_iw>w.
\end{cases}$. But $\partial_i c_w(x)$ is fixed by $s_i$, so
$$\partial_i c_{w}(x)=\begin{cases}
c_{s_iw}(x), & s_iw<w,\\
0,& s_iw>w.
\end{cases}. $$
It is clear that $c_{w_0}(x)=F(x)$.
So the lemma follows from induction.

\paragraph{Proof of Theorem \ref{doubleidentities}. }
Now, we are going to prove two identities in Theorem \ref{doubleidentities}.
We apply Lemma \ref{TotalLeibniz} on
$\mathfrak{S}_{w_0}(w_0y,x)\mathfrak{S}_{w}(y,t)$,
$$\begin{array}{rl}
 (-1)^{w_0}\partial_{w_0}^y\big(\mathfrak{S}_{w_0}(w_0y,x)\mathfrak{S}_{w}(y,t)\big)
&\equiv \displaystyle \sum_{u}(-1)^u\big(\partial_{uw_0}^y\mathfrak{S}_{w_0}(y,x)\big)\cdot \partial_u^y\mathfrak{S}_{w}(y,t)\\
&\equiv \displaystyle \sum_{w=wu^{-1}\odot u}
(-1)^u\mathfrak{S}_{u^{-1}}(y,x)\cdot \mathfrak{S}_{wu^{-1}}(y,t)\\
&\equiv \displaystyle \sum_{w=u\odot v}
(-1)^v\mathfrak{S}_{v^{-1}}(y,x)\cdot \mathfrak{S}_{u}(y,t).
\end{array}$$
Since after applying $\partial_{w_0}$, any polynomial becomes symmetric, so the right hand side is symmetric in $y$.
But the ideal factored allow us to write $f(x,y,t)=f(x,x,t)=f(x,t,t)$ for polynomial $f$ symmetric in $y$. So in particular,
$$\begin{array}{l}
\quad \displaystyle \sum_{w=u\odot v}
(-1)^v\mathfrak{S}_{v^{-1}}(y,x)\cdot \mathfrak{S}_{u}(y,t)\\
\equiv \displaystyle \sum_{w=u\odot v}
(-1)^v\mathfrak{S}_{v^{-1}}(x,x)\cdot \mathfrak{S}_{u}(x,t)\equiv \mathfrak{S}_{w}(x,t),\\
\equiv \displaystyle \sum_{w=u\odot v}
(-1)^v\mathfrak{S}_{v^{-1}}(t,x)\cdot \mathfrak{S}_{u}(t,t)\equiv (-1)^w\mathfrak{S}_{w}(t,x).\\
\end{array}$$
Since $\mathfrak{S}_w(x,x)=\begin{cases}1,& w=1\\ 0, & w\neq 1.\end{cases}$.
This is the second identity. Taking in the second identity to above, finally we get
$$\mathfrak{S}_{w}(x,t)\equiv \sum_{w=u\odot v}\mathfrak{S}_{v}(x,y)\cdot \mathfrak{S}_{u}(y,t). $$
The proof is complete.

\section{Relation to Convolution Algebra}

In this section, I will give a geometric proof of Lemma \ref{TotalLeibniz}.

\paragraph{Convolution Algebra. }
Lemma \ref{TotalLeibniz} can also be proved taking advantage of the convolution in geometric representation theory.
It is generally more convenient to use equivariant Borel--Moore homology see \cite{chriss2009representation},
rather than equivariant cohomology (while they are dual).
In our specific situation, it is useful to use cohomology.
Let $p_{i}: G/B\times G/B\times G/B\to G/B\times G/B$ be the projection by omitting $i$-factor.
The convolution operator $*$ over $H^*_G(G/B\times G/B)$ is defined by the following diagram
$$\xymatrix{
\mbox{$\begin{array}{l}
H_G^*(G/B\times G/B)\\
\quad \otimes H_G^*(G/B\times G/B)
\end{array}$}\ar[rr]^{p_3^*\otimes p_1^*}\ar[d]_{*}&&
\mbox{$\begin{array}{l}
H_G^*(G/B\times G/B\times G/B)\\
\quad \otimes H_G^*(G/B\times G/B\times G/B)
\end{array}$}\ar[d]^{\smile}\\
H_G^*(G/B\times G/B)&&H_G^*(G/B\times G/B\times G/B)\ar[ll]^{(p_2)_*}.
}$$

Note that
$$\begin{array}{rl}
H_G^*(G/B\times G/B)&=H^*(\pt\underset{G}* (G/B\times G/B))\\
&=H^*(\pt\underset{G}* G\times_B G/B)\\
&=H^*(\pt\underset{B}* G/B)\\
&=H^*(\pt\underset{B}*G\underset{B}*\pt);\\
H_G^*(G/B\times G/B\times G/B)&=H^*(\pt\underset{G}* (G/B\times G/B\times G/B))\\
&=H^*(\pt\underset{G}* G\times_B G\times_B G/B)\\
&=H^*(\pt\underset{B}*G\underset{B}*G\underset{B}*\pt)\\
\end{array}$$
In this case, $p_3$ and $p_1$ are exactly $\pi_1,\pi_2$ defined before, so
$$\begin{array}{rl}
p_3^*& H^*(\pt\underset{B}*G\underset{B}*\pt) \to H^*(\pt\underset{B}*G\underset{B}*G\underset{B}*\pt)
\qquad f(x,t)\longmapsto f(y,t)\\
p_1^*& H^*(\pt\underset{B}*G\underset{B}*\pt) \to H^*(\pt\underset{B}*G\underset{B}*G\underset{B}*\pt)
\qquad f(x,t)\longmapsto f(x,y)\\
\end{array}$$
From the map and the functoriality of push forward
$$\xymatrix{
(G/B\times G/B\times G/B)_G\ar[d]\ar[r]& \pt\underset{B}* G\underset{B}* G\underset{B}*\pt \ar[d]\ar[r]&
\pt\underset{G}* G\underset{B}* G\underset{G}*\pt\ar[d]\\
(G/B\times G/G \times G/B)_G\ar[r]& \pt\underset{B}*G \underset{G}* G\underset{B}*\pt\ar[r]&
 \pt\underset{G}*G \underset{G}* G\underset{G}*\pt, \\
}$$
we see the push forward takes the form
$$(p_2)_*: H^*(\pt\underset{B}*G\underset{B}*G\underset{B}*\pt) \to H^*(\pt\underset{B}*G\underset{B}*\pt) \qquad
f(x,y,t)\longmapsto \partial^y_{w_0}f(x,y,t)\big|_{y=t}. $$
So the convolution takes the form
$$f*g=\partial_{w_0}^y g(x,y)f(y,t)\big|_{y=t}. $$
Note that we can define the action of $H_G^*(G/B\times G/B)$ on $H^*_G(G/B)=H^*(BT)\cong \mathbb{Z}[t_1,\ldots,t_n]$
which makes $H_G^*(G/B)$ an $H_G^*(G/B\times G/B)$-$H_G^*(\pt)$-bimodule.
In this case the convolution has the same expression.
$$f*g=\partial_{w_0}^y g(y)f(y,t)\big|_{y=t}. $$
Since the action is faithful from the expression (it does depend on $x$), so
Lemma \ref{TotalLeibniz} is equivalent to
$$f*g=\displaystyle (-1)^{w_0}\sum_{u\in W} (-1)^u \partial_{u}^yg(y)\partial_{uw_0}^yf((w_0y,t))\big|_{y=t}$$
for all $f\in H_B^*(G/B)$ and $g\in H_G(G/B)$.

\paragraph{Schubert Cells in $H_G(G/B\times G/B)$. }
Since
$$\pt \underset{B}*G/B\times G/B=\pt\underset{B}*G/B, $$
the $G$-orbits of $G/B\times G/B$ are
in bijection with $B$-orbits of $G/B$.
For $w\in W$,
denote $\Psi_w=BwB/B$ the Schubert cell,
and $\Lambda_w$ the $G$-orbit correspondent to $BwB/B$, explicitly,
$$\Lambda_w=\{(xB,yB): x^{-1}y\in BwB\}. $$

Let $[\Psi_w]_B$ the equivariant cohomology class of $\overline{\Psi_w}$ in $G/B$,
and $[\Lambda_w]_G$ the equivariant cohomology class of $\overline{\Lambda_w}$ in $G/B\times G/B$.
They are the same under the isomorphism $H^*_B(G/B)=H_G^*(G/B\times G/B)$ from the definition of fundamental class.

On one hand, it is well-known \cite{ginzburg1998geometric}
that $[\Lambda_w]_G$ acts as the Demazure operator $\partial_w$ over $H_G^*(G/B)$.
On the other hand, $[\Psi_w]_B=\mathfrak{S}_{w_0w}(x,w_0t)$ from our computation.
In particular, $[\Psi_w]_B$ forms a basis of $H_B(G/B;\mathbb{Q})$ as $H_G(G/B;\mathbb{Q})$-module.
So it suffices to check the expression of convolution above for the case $f=\mathfrak{S}_{w_0w}(x,w_0t)$.
Besides, it is easy to check the following commuting diagram
$$\xymatrix{
\mathbb{Z}[t]\ar[r]^{w_0}\ar[d]_{(-1)^v\partial_v}
& \mathbb{Z}[t]\ar[d]^{\partial_{w_0vw_0}}\\
\mathbb{Z}[t]\ar[r]_{w_0}& \mathbb{Z}[t]
}$$
So,
$$\begin{array}{l}
\displaystyle\quad (-1)^{w_0}\sum_{u\in W} (-1)^u \partial_{u}^yg(y)\partial_{uw_0}^y(\mathfrak{S}_{w_0w}(w_0y,w_0t))\big|_{y=t}\\
\displaystyle= (-1)^{w_0}\sum_{u\in W} (-1)^u \partial_{u}^yg(y)(-1)^{w_0u}(\partial_{w_0u}^y\mathfrak{S}_{w_0w})((w_0y,w_0t))\big|_{y=t}\\
= \partial_{u}^yg(y).
\end{array}$$
This is the geometric proof of Lemma \ref{TotalLeibniz}.

\bibliographystyle{plain}
\bibliography{bibfile}

\end{document}